
\magnification=\magstep1
\def\refer#1{ \medskip \parindent=34pt \hang \noindent \rlap{[#1]\hfil} \hskip 30pt }
\baselineskip=14pt

\font\tenmsb=msbm10
\font\sevenmsb=msbm7
\font\fivemsb=msbm5
\newfam\msbfam
\textfont\msbfam=\tenmsb
\scriptfont\msbfam=\sevenmsb
\scriptscriptfont\msbfam=\fivemsb
\def\Bbb#1{{\fam\msbfam\relax#1}}
\def\lam{{\lambda}}

\def\al{{\alpha}}

\def\proof{{\medskip\noindent {\bf Proof. }}}

\def\qed{{\hfill $\square$ \bigskip}}
\def\subsec#1{{\bigskip\noindent \bf{#1}.}}

\def\cite#1{{[#1]}}

\def\eps{\varepsilon}

\def\angel#1{{\langle#1\rangle}}
\def\norm#1{\Vert #1 \Vert}

\def\q {\quad} \def\qq {\qquad}
\def\frac#1#2{{#1\over #2}}
\def\del{{\partial}}
\def\wt{\widetilde}

\def\wh{\widehat}

\def\ni{\noindent }
\def\ms{\medskip}
\def\bs{\bigskip}
\def\cl#1{\centerline{#1}}

\parindent=30pt

\def\Cov{{\mathop {{\rm Cov\, }}}}

\def\square{{\vcenter{\vbox{\hrule height.3pt
\hbox{\vrule width.3pt height5pt \kern5pt
Ê\vrule width.3pt}
\hrule height.3pt}}}}

\def\tfrac#1#2{{\textstyle {#1\over #2}}}

\def\tlint{{- \kern-0.85em \int \kern-0.2em}} 
\def\dlint{{- \kern-1.05em \int \kern-0.4em}}

\def\y0{\langle y_0,\be_i\rangle}
\def\R{{\Bbb R}}
\def\E{{{\Bbb E}\,}}
\def\P{{\Bbb P}}

\def\F{{\cal F}}

 \def\sB {{\cal B}} 
\def\sD {{\cal D}}  
  
  \def\sL {{\cal L}}

\def\sS {{\cal S}} \def\sT {{\cal T}}

 \def\bN {{\Bbb N}}

\def\be{{\epsilon}}
\def\bal{{S^\alpha}}

\def\bpto{\smash{
Ê\mathop{\longrightarrow}\limits^{bp}}}

\font\ften=cmr8

\cl{\bf Infinite dimensional stochastic differential equations}
\ms
\cl{\bf of Ornstein-Uhlenbeck type}

\bs

\centerline{Siva R. Athreya,\q Richard F. Bass${}^{1}$, 
\q Maria Gordina${}^{2}$
\q and \q Edwin A. Perkins${}^{3}$}

\footnote{$\empty$}{\ften 1. Research supported in part by NSF grant DMS0244737.}

\footnote{$\empty$}{\ften 2. Research supported in part by ÊNSF grant DMS0306468.}
\footnote{$\empty$}{\ften 3. Research supported in part by ÊNSERC Research Grant.}
\footnote{$\empty$}{\ften {Keywords:} Semigroups, H\"older spaces, perturbations, resolvents, }
\footnote{$\empty$}{\ften \q \q \q \q \q elliptic operators, Ornstein-Uhlenbeck processes,}
\footnote{$\empty$}{\ften \q \q \q \q\q infinite dimensional stochastic differential equations}
\footnote{$\empty$}{\ften Classification: Primary 60H10. Secondary 35R15, 60H30}
\footnote{$\empty$}{\ften Running Head: SDEs of O-U type}

\vskip0.5truein

{ \centerline{\bf Abstract}} 
We consider the operator
$$\sL f(x)=\tfrac12 \sum_{i,j=1}^\infty a_{ij}(x)\frac{\del^2 f}{\del x_i \del x_j}(x)-\sum_{i=1}^\infty \lam_i x_i b_i(x) \frac{\del f}{\del x_i}(x).$$
We prove existence and uniqueness of solutions to the martingale problem
for this operator under appropriate conditions
on the $a_{ij}, b_i$, and $\lam_i$.
The process corresponding to $\sL$ solves an infinite dimensional stochastic
differential equation similar to that for the infinite dimensional Ornstein-Uhlenbeck
process.

\vfil\eject

\subsec{1. Introduction}

Let $\lam_i$ be a sequence of positive reals tending to infinity, let $\sigma_{ij}$
and $b_i$ be functions defined on a suitable Hilbert space
which satisfy certain continuity and non-degeneracy conditions,
and let $W^i_t$ be a sequence of independent one-dimensional Brownian motions. In
this paper we consider the countable system of stochastic differential equations
$$dX^i_t=\sum_{j=1}^\infty \sigma_{ij}(X_t) dW^i_t-\lam_i b_i(X_t) X^i_t\, dt, 
\qq i=1,2,\ldots, \eqno (1.1)$$
and investigate sufficient conditions for weak existence and weak uniqueness to hold.
Note that when the $\sigma_{ij}$ and $b_i$ are constant, we have the stochastic
differential equations characterizing the infinite-dimensional Ornstein-Uhlenbeck process.

We approach the weak existence and uniqueness of (1.1) by means of the martingale 
problem for the corresponding operator
$$\sL f(x)=\tfrac12 \sum_{i,j=1}^\infty a_{ij}(x) \frac{\del^2f}{\del x_i \del x_j}(x)-
\sum_{i=1}^\infty \lam_i x_i b_i(x) \frac{\del f}{\del x_i}(x) \eqno (1.2)$$
operating on a suitable class of functions, where $a_{ij}(x)=\sum_{k=1}^\infty \sigma_{ik}(x)
\sigma_{jk}(x)$. Our main theorem says that if the 
$a_{ij}$ are nondegenerate and bounded, the $b_i$ are bounded above and below, 
Êand the $a_{ij}$ and $b_i$ satisfy
appropriate H\"older continuity conditions, then existence and uniqueness hold for
the martingale problem for $\sL$; see Theorem 5.7 for a precise statement.

There has been considerable interest in infinite dimensional operators
whose coefficients are only H\"older continuous. For perturbations of the Laplacian,
see Cannarsa and DaPrato [CD], where Schauder estimates are proved using interpolation
theory and then applied to Poisson's equation in infinite dimensions with H\"older
continuous coefficients (see also [DZ]). 

Similar techniques have been used to study
operators of the form (1.2). In finite dimensions see Ê
[L1], [L2], [L3], and [DL]. For the infinite dimensional case see [Ce], [D1], [D2], [D3], 
[D4], [DZ], and [Z]. Common to all of these papers is the use of interpolation theory to
obtain the necessary Schauder estimates. In functional analytic terms,
the system of equations (1.1) is a special case of the equation
$$dX=(b(X_t) X_t+F(X_t)) \, dt +\sqrt {a(X_t)}\, dW_t, \eqno (1.3)$$
where $a$ Ê
is a mapping from a Hilbert space $H$ to the space of
bounded non-negative self-adjoint linear operators on $H$, 
$b$ is a mapping from $H$ to the non-negative self-adjoint linear operators on $H$
(not necessarily bounded),
$F$ is a bounded operator on $H$,
and
$b(x)x$ represents the composition of operators.
Previous work has concentrated on (1.3) in the following cases: where $a$ is constant, Ê$b$ is Lipschitz continuous,
and $F\equiv 0$;
where $a$ and $b$ are constant and $F$ is bounded; 
and Ê where $F$ is bounded, $b$ is constant and $a$ is a Êperturbation
of a constant operator by means of a H\"older continuous nonnegative self-adjoint operator.
We also mention the paper [DR] Êwhere weak solutions to (1.3) are considered.
In our paper we consider the
equation (1.3) with the $a$ and $b$ satisfying certain H\"older conditions and $F\equiv 0$. There
would be no difficulty introducing bounded $F(X_t) dt$ terms, but we chose
not to do so.

The paper most closely related to this one is that of Zambotti [Z].
Our results complement those of [Z] as each has its own advantages.
We were able to remove the restriction that the $a_{ij}$'s be given by means
of a perturbation by a bounded nonnegative operator which in turn facilitates localization, but at
the expense of working with respect to a fixed basis and hence imposing summability conditions
involving the off-diagonal $a_{ij}$. 
See Remark 5.10 for a further discussion in light of a couple of examples and our explicit
hypotheses for Theorem~5.7.

There are also martingale problems for infinite dimensional operators with 
H\"older continuous coefficients that arise from
the fields of superprocesses and stochastic partial differential equations (SPDE). 
See [P] for a detailed introduction to these. We mention [DM], where superprocesses
in the Fleming-Viot setting are considered, and [BP], where uniqueness of a martingale
problem for superprocesses on
countable Markov chains with interactive branching is shown to hold. 
These latter results motivated the present approach as the weighted
H\"older spaces used there for our perturbation bounds coincide with the function spaces 
$S^\alpha$ used here (see Section 2), at least in the finite-dimensional setting (see [ABP]).

Consider the one dimensional SPDE
$$\frac{\del u}{\del t}(t,x)=\tfrac12 \frac{\del^2 u}{\del x^2}(x,t) +A(u) d\dot W, \eqno (1.5)$$
where $\dot W$ is space-time white noise. If one sets
$$X^j_t=\int_0^{2\pi} e^{ijx} u(x,t)\, dx, \qq j=0, \pm1, \pm2, \ldots,$$
then the collection $\{X^i\}_{i=-\infty}^\infty$ can be shown to solve Êthe system
(1.1) with $\lam_i= i^2$, the $b_i$ constant, and the $a_{ij}$ defined
in an explicit way in terms of $A$. Our original interest in the problem solved
in this paper was to understand (1.5) when the coefficients $A$ were bounded above and below
but were only H\"older continuous as a function of $u$. The results in this paper do not
apply to (1.5) and we hope to return to this in the future. Ê

\ms
The main novelties of our paper are the following.

\item{(1)} {\sl $C^\al$ estimates (i.e., Schauder estimates) for the infinite dimensional
Ornstein-Uhl\-en\-beck process.} \quad
These were already known Ê(see [DZ]), but we point out that in contrast to using interpolation
theory, our derivation is quite elementary and relies on a simple real variable lemma
together with some semigroup manipulations. 

\item{(2)} {\sl Localization.} \quad
We use perturbation
theory Êalong the lines of Stroock-Varadhan to establish uniqueness of the martingale
problem when the coefficients are sufficiently close to constant. We then perform
a localization procedure to establish our main result. In infinite dimensions localization
is much more involved, and this argument represents an important feature of this work.

\item{(3)} {\sl A larger class of perturbations.} \quad Unlike much of the previous work
cited above,
we do not require that the perturbation of the second order term be by an operator that
is nonnegative. The price we pay is that we require additional conditions
on the off-diagonal $a_{ij}$'s.
\ms

After some definitions and preliminaries in Section 2, we establish the needed Schauder
estimates in Section 3. Section 4 contains the proof of existence and Section 5 the
uniqueness. Section 5 also contains some specific examples where our main result applies. Ê
This includes coefficients $a_{ij}$ which depend on a finite number of local coordinates near
$(i,j)$ in a H\"older manner.

We use the letter $c$ with or without subscripts for
finite positive constants whose value is unimportant and which may
vary from proposition to proposition. $\alpha$ will denote a real number
between $0$ and $1$.

\ms
\ni{\bf Acknowledgment.} We would like to thank L. Zambotti for helpful
conversations on the subject of this paper.

\subsec{2. Preliminaries}
We use the following notation. ÊIf $H$ is a separable Hilbert space
and $f:H\to \R$, $D_wf(x)$ is the directional derivative of $f$ at
$x\in H$ in the direction $w$; we do not require $w$ to be a unit
vector. ÊThe inner product in $H$ is denoted $\angel{\cdot, \cdot}$,
and $|\cdot |$ denotes the norm generated by this inner
product. $C_b=C_b(H)$ is the collection of $\R$-valued bounded
continuous functions on $H$ with the usual supremum norm.
Let $C_b^2$ be the set of
functions in $C_b$ for which the first and second order partials are
also in $C_b$. For $\al\in (0,1)$, set
$$|f|_{C^\al}=\sup_{x\in H, h\ne 0} \frac{|f(x+h)-f(x)|}{|h|^\al}$$
and let $C^\al$ be the set of functions in $C_b$ for which $\Vert
f\Vert_{C^\al}=\norm{f}_{C_b}+|f|_{C^al}$
is finite. 

Let $V:\sD(V)\to H$ be a (densely defined)
self-adjoint non-negative definite operator such that
$$V^{-1}\hbox{ is a trace class operator on }H,\eqno (2.1)$$
Then there is a complete orthonormal system $\{\be_n:n\in \bN\}$ of eigenvectors of
$V^{-1}$ with corresponding eigenvalues $\lambda_n^{-1}$, $\lambda_n>
0$, satisfying 
$$\sum_{n=1}^\infty \lambda_n^{-1}<\infty, \q \lambda_n\uparrow\infty,
\q V\be_n=\lambda_n\be_n$$ (see, e.g. Section 120 in [RN]). Let $Q_t=e^{-tV}$ be the
semigroup of contraction operators on $H$ with generator $-V$. If $w\in H$, 
let $w_n=\langle w, \be_n\rangle$ and we
will write $D_if$ and $D_{ij}f$ for $D_{\be_i}f$ and $D_{\be_i}D_{\be_j}f$,
respectively.

Assume $a:H\to L(H,H)$ is a mapping from $H$ to the space of
bounded self-adjoint operators on $H$ and $b:H\to L(\sD(V),H)$ is
a mapping from $H$ to self-adjoint non-negative definite operators
on $\sD(V)$ such that $\{\be_n\}$ are eigenvectors of $b(x)$ for
all $x\in H$. ÊIf 
$a_{ij}(x)=\langle \be_i,a(x)\be_j\rangle$ and
$b(x)(\be_i)=\lambda_ib_i(x)\be_i$, we assume that for
some
$\gamma>0$
$$\eqalign{\gamma^{-1} |z|^2&\geq \sum_{i,j} 
a_{ij}(x)z_iz_j\geq \gamma |z|^2, 
\qq x,z\in H,\cr
\gamma^{-1}&\geq b_i(x)\geq \gamma, \qq x\in H, \quad
i\in\bN.\cr}\eqno (2.2)$$

We consider the martingale problem Êfor the operator $\sL$ which, with
respect to the coordinates $\langle x,\be_i\rangle$, is 
defined by 
$$\sL f(x)=\tfrac12 \sum_{i,j=1}^\infty a_{ij}(x)D_{ij}
f(x)-\sum_{i=1}^\infty \lam_{i} x_i b_i(x)D_i f(x). \eqno (2.3)$$
Let $\sT$ be the class of functions in $C_b^2$ that depend on only finitely
many coordinates and $\sT_0$ be the set of functions in $\sT$ with compact support.
More precisely, $f\in \sT$ if there exists $n$ and $f_n\in C_b^2(\R^n)$
such that $f(x_1, \ldots, x_n, \ldots)=f_n(x_1, \ldots x_n)$
for each point $(x_1, x_2, \ldots)$ and $f\in\sT_0$ if, in addition, $f_n$ has
compact support. Let $X_t$ denote the coordinate maps on the space 
$C([0,\infty),H)$ of continuous
$H$-valued paths. ÊWe say that a
probability measure $\P$ on $C([0,\infty),H)$ is a solution to the martingale
problem for $\sL$ started at
$x_0$ if 
$\P(X_0=x_0)=1$ and 
$f(X_t)-f(X_0)-\int_0^t \sL f(X_s) ds$ is a martingale for
each $f\in \sT$.

The connection between systems of stochastic differential equations
and martingale problems continues to hold in infinite dimensions; see,
for example,
[KX] pp.~166--168. We will use this fact without further mention.

There are different possible martingale problems depending
on what class of functions we choose as test functions. Since
existence is the easier part 
for the martingale problem (see Theorem 4.2 below)
and uniqueness is the more difficult part, we will get a
stronger and more useful theorem if we have a smaller
class of test functions. The collection $\sT$ is a reasonably small class. 
When $a(x)\equiv a^0$ and $b(x)\equiv V$ are constant functions,
the process associated with $\sL$ is the well-known
$H$-valued Ornstein-Uhlenbeck process. We briefly recall the definition;
see Section 5 of [ABP] for details. Let $(W_t,t\ge 0)$ be the
cylindrical Brownian motion on $H$ with covariance
$a$. 
Let $\F_t$ be the right continuous filtration generated by $W$.
Consider the stochastic differential equation
$$dX_t=dW_t-VX_t\,dt.\eqno (2.4)$$
There is a pathwise unique solution to (2.4) 
whose laws
$\{\P^x,x\in H\}$ define a unique homogeneous strong Markov process on the
space of continuous
$H$-valued paths (see, e.g. Section 5.2 of [KX]). 
$\{X_t,t\ge 0\}$ is an $H$-valued Gaussian process satisfying
$$\E(\langle X_t,h\rangle)=\langle X_0,Q_th\rangle \hbox{ for all }
h\in H,\eqno (2.5)$$ and 
$$\Cov(\langle X_t,g\rangle\langle X_t,h\rangle)=\int_0^t\langle
Q_{t-s}h,aQ_{t-s}g\rangle ds.\eqno(2.6)$$
The law of $X$ Êstarted at $x$ solves the 
martingale problem for
$$\sL_0 f(x)=\tfrac12 \sum_{i,j=1}^\infty a_{ij}^0D_{ij}
f(x)-\sum_{i=1}^\infty \lam_{i} x_i D_i f(x). \eqno (2.7)$$

We let $P_tf(x)=\E^x f(X_t)$ be the semigroup corresponding to $\sL_0$,
and \hfil\break
$R_\lam=\int_0^\infty e^{-\lam s} P_s\, ds$ be the
corresponding resolvent. We define the semigroup norm 
Ê$\norm{\cdot}_\bal$ for $\al\in (0,1)$ by
$$|{f}|_\bal= \sup_{t>0} t^{-\al/2}
\norm{P_tf-f}_{C_b}\eqno (2.8)$$
and
$$\norm{f}_\bal=\norm{f}_{C_b} +|f|_\bal.$$
Let $S^\al$ denote the space of measurable functions on $H$ for which
this norm is finite. 

For $x\in H$ and $\beta\in(0,1)$ define $|x|_{\beta}=\sup_k|\langle x,
\be_k\rangle|\lam_k^{\beta/2}$ and $$H_{\beta}=\{x\in H:|x|_{\beta}<\infty\}. \eqno (2.9)$$

\subsec{3. Estimates}

We start with the following real variable lemma.

\proclaim Lemma 3.1. Let $A >0, B>0.$ Assume $K: C_b(H)\to C_b(H)$ is a bounded linear operator such that $$\norm{Kf}_{C_b} \leq A\norm{f}_{C_b},
\qq f\in {C_b}(H),\eqno (3.1)$$ and there exists $v\in H$ such that
$$\norm{Kf}_{C_b}\leq B\norm{D_vf}_{C_b}, \eqno (3.2)$$ for all $f$
such that $D_vf\in C_b(H)$. Then for each $\alpha\in (0,1)$ there is a constant $c_1=c_1(\alpha)$ such that
$$\norm{Kf}_{C_b}\leq c_1 |v|^\al |f|_{C^\al} B^\al A^{1-\al} \qq \hbox{for all }f\in
C^\al.$$

\proof Assume (3.1) and (3.2), the latter for some $v\in H$. Let $\{p_t:
t\geq 0\}$ be the standard Brownian density on
$\R$. If $f\in C^\al$, set
$$p_\eps*f(x)=\int_\R f(x+zv) p_\eps(z)\, dz, \qq x\in H.$$
Since Êa change of variables shows that
$$p_\eps*f(x+hv)-p_\eps*f(x)=\int_\R f(x+zv)p_\eps(z-h)\, dz-\int_\R f(x+zv)p_\eps(z)\, dz,$$
it follows that 
$$D_v(p_\eps*f)(x)=-\int f(x+zv) p'_\eps(z)\, dz;$$
this is in $C_b(H)$
and
$$\eqalign{|D_v(p_\eps*f)(x)|&=\Bigl|-\int f(x+zv)p_\eps'(z)\, dz\Bigr|\cr
&=\Bigl|\int (f(x+zv)-f(x))p_\eps'(z)\, dz\Bigr|\cr
&\leq |f|_{C^\al} |v|^\al \int |z|^\al \frac{|z|}{\eps} p_\eps (z)\, dz
\cr
&=c_2|f|_{C^\al} |v|^\al \eps^{(\al-1)/2},\cr}$$
where
$c_2=\int |z|^{\al+1} p_1(z)\, dz$. We therefore obtain from (3.2) that
$$\norm{K(p_\eps*f)}_{C_b}\leq c_2B|f|_{C^\al} |v|^\al \eps^{(\al-1)/2}.\eqno (3.3)$$
Next note that
$$\eqalign{|p_\eps*f(x)-f(x)|&\leq \int |f(x+zv)-f(x)|p_\eps(z)\, dz\cr
&\leq |f|_{C^\al} |v|^\al \int|z|^\al p_\eps (z) \, dz\cr
&=c_3|f|_{C^\al} |v|^\al \eps^{\al/2},\cr}$$
where $c_3=\int |z|^\al p_1(z)\, dz$. By (3.1) 
$$\norm{K(p_\eps *f-f)}_{C_b}\leq c_3A |f|_{C^\al} |v|^\al \eps^{\al/2}. \eqno (3.4)$$
Let $c_4=c_2\lor c_3$ and $\eps=B^2/A^2$. Combining (3.3) and (3.4)
we have
$$\eqalign{\norm{Kf}_{C_b}&\leq c_4|f|_{C^\al} |v|^\al \eps^{\al/2} [A+B\eps^{-1/2}]\cr
&=2c_4|f|_{C^\al} |v|^\al B^\al A^{1-\al}.\cr}$$
\qed

Set $$h(u)=\cases{(2u)/(e^{2u}-1) & $u\ne 0$;\cr
1 & $u=0$.}$$
and 
$$|w|_t=\Big(\sum_i w_i^2 h(\lam_i t)\Big)^{1/2} \leq |w|.$$
Recall $$Q_t w=\sum_{i=1}^{\infty} e^{-\lam_i t} w_i e_i.$$
We have the following by Propositions 5.1 and Ê5.2 of [ABP]:

\proclaim Proposition 3.2. (a) For all $w\in H$, $f\in {C_b}(H)$, and $t>0$,
$D_wP_tf\in C_b(H)$ and
$$\norm{D_wP_t f}_{C_b} \leq \frac{|w|_t \norm{f}_{C_b}}{\sqrt{\gamma t}}. \eqno (3.5)$$
\item{(b)} If $t\geq 0$, $w\in H$, and $f:H\to \R$ is in ${C_b}(H)$ such that
$D_{Q_tw}f\in C_b(H)$, then
$$D_wP_tf(x)=P_t(D_{Q_tw}f)(x), \qq x\in H.$$
In particular,
$$\norm{D_wP_tf}_{C_b}\leq \norm{D_{Q_tw}f}_{C_b}.\eqno (3.6)$$

We now prove

\proclaim Corollary 3.3. Let $f\in C^\al$, $u,w\in H$. Then 
for all $t>0$, $D_wP_tf$ and $ D_uD_wP_tf$ are in $C_b(H)$
and there exists a constant $c_1=c_1(\alpha,\gamma)$ independent of $t$ such that
$$\norm{D_wP_tf}_{C_b}\leq c_1 |w|_t|f|_{C^\al} t^{(\al-1)/2}\leq
c_1|w|\, |f|_{C^\al} t^{(\al-1)/2}, \eqno (3.7)$$
and 
$$\eqalignno{ \norm{D_uD_wP_tf}_{C_b}&\leq c_1 |Q_{t/2}u|_{t/2} |w|_{t/2}
|f|_{C^\al} t^{\frac{\al}{2}-1}\leq c_1|u|_{t/2}|w|_{t/2} |f|_{C^\al} t^{\frac{\al}{2}-1}
\cr
&\leq c_1|u|\, |w|\, |f|_{C^\al} t^{\frac{\al}{2}-1}. &(3.8)\cr}$$

\proof That $D_wP_tf$ is in $C_b(H)$ is immediate from Proposition~3.2(a).
By (3.5) and (3.6) we may
Êapply
Lemma 3.1 to $K=D_wP_t$ with $v=Q_tw$, $A=|w|_t(\gamma t)^{-1/2}$ and $B=1$ 
to conclude for $f\in C^\al$
$$\eqalignno{\norm{D_wP_tf}_{C_b}&\leq c_2|Q_tw|^\al |f|_{C^\al} 
|w|_t^{1-\al} (\gamma t)^{-(1-\al)/2}\cr
&\leq c_2 \gamma^{(\al-1)/2}|w|_t|f|_{C^\al} t^{(\al-1)/2}.&(3.9)\cr}$$
This gives (3.7).

By Proposition 3.2, $D_wD_uP_tf=D_wP_{t/2}D_{Q_{t/2}u}P_{t/2}f$, and the
latter is seen to be in $C_b(H)$ by invoking Proposition~3.2(a)
twice. ÊUsing (3.5) and then (3.9) we have
$$\eqalign{\norm{D_wD_uP_tf}_{C_b}&=\norm{D_wP_{t/2}D_{Q_{t/2}u}P_{t/2}f}_{C_b}\cr
&\leq |w|_{t/2} (\gamma t/2)^{-1/2} \norm{D_{Q_{t/2}u}P_{t/2}f}_{C_b} 
\cr
&\leq |w|_{t/2}(\gamma t/2)^{-1/2} c_2 \gamma^{(\al-1)/2} |Q_{t/2}u|_{t/2}
|f|_{C^\al} ( t/2)^{(\al-1)/2}.\cr}$$
This gives (3.8). \qed

\noindent {\bf Remark 3.4.} We often will use the fact that
there exists $c_1$ such that
$$\norm{f}_{C^\al}\leq c_1 \norm{f}_\bal. \eqno (3.10)$$
This is (5.20) of [ABP].

\proclaim Corollary 3.5. There exists $c_1=c_1(\al, \gamma)$ such that
for all $\lam>0$, $f\in C^{\alpha}$, $i\leq j$,
we have $D_iR_\lam f, D_{ij}R_\lam f\in C_b$, and
$$\norm{D_iR_\lam f}_{C_b}\leq c_1(\lam+\lam_i)^{-(\al+1)/2}|f|_{C^\al}.
\eqno (3.11)$$
$$\norm{D_{ij}R_\lam f}_{C_b} \leq c_1(\lam+\lam_j)^{-\al/2}
|f|_{C^\al}.\eqno (3.12)$$
$$\norm{D_{i}R_\lam f}_{C^\al}
\leq c_1(\lam+\lam_i)^{-1/2} \norm{f}_{C^\al}. \eqno (3.13)$$
$$\norm{D_{ij}R_\lam f}_{C^\al}
\leq c_1 \norm{f}_{C^\al}. \eqno (3.14)$$

\proof 
Corollary 3.3 is exactly the same as Proposition 5.4 in [ABP], but
with the $S^\al$ norms replaced by $C^\al$ norms. We may therefore 
follow the proofs of Theorem 5.6 and Corollary 5.7 in [ABP] and then
use (3.10) to obtain our result. 
However the proofs in [ABP] can be streamlined, so for the sake of
clarity and completeness we give a more straightforward proof. 

From (3.7) and (3.8) we may differentiate under the time integral and
conclude that the first and second order partial derivatives of
$R_\lambda f$ are continuous. To derive (3.12), note first that by (3.8), 
$$\eqalignno{\norm{D_{ij}P_t f}_{C_b}&=\norm{D_{ji}P_tf}_{C_b} \leq c_2|Q_{t/2}
\be_j|\, |\be_i| \, |f|_{C^\al} t^{\frac{\al}{2}-1}&(3.15)\cr
&= c_2 e^{-\lam_j t/2} |f|_{C^\al} t^{\frac{\al}{2}-1}.\cr}$$
Multiplying by $e^{-\lam t}$ and integrating over $t$ from 0 to $\infty$
yields (3.12).

Next we turn to (3.14). Recall the definition of the $S^\al$ norm from
(2.8). In view Êof (3.10) it suffices to show
$$\norm{D_{ij} R_\lam f}_{S^\al} \leq c_3\norm{f}_{C^\al}.$$
Since 
$$\norm{P_tD_{ij}R_\lam f-D_{ij}R_\lam f}_{C_b}\leq 2\norm{D_{ij}R_\lam f}_{C_b}
\leq c_1 |f|_{C^\al} (\lam +\lam_j)^{-\al/2}$$
by (3.12), we need only consider $t\leq (\lam+\lam_j)^{-1}$.

Use Proposition 3.2(b) to write
$$\eqalignno{P_tD_{ij} R_\lam f&-D_{ij}R_\lam f =[e^{-\lam_i t} e^{-\lam_j t} D_{ij} P_t R_\lam f
-D_{ij} P_t R_\lam f]\cr
&\qq + [D_{ij}P_t R_\lam f-D_{ij} R_\lam f].&(3.16)\cr}$$
Recalling that $\lambda_i\le \lambda_j$, we see that 
the first term is bounded in absolute value by
$$\eqalign{c_4(\lam_jt)^{\al/2} \norm{D_{ij} P_t R_\lam f}_{C_b}
&\leq c_5t^{\al/2} \int_0^\infty \lam_j^{\al/2} e^{-\lam s} \norm{D_{ij} P_{t+s}f}_{C_b}
ds\cr
&\leq c_5t^{\al/2} |f|_{C^\al},\cr}$$
using (3.15).

The second term in (3.16) is equal, by the semigroup property, to
$$\eqalign{\int_0^\infty e^{-\lam s} D_{ij} 
&P_{t+s} f \, ds-\int_0^\infty e^{-\lam s} D_{ij} P_s f\, ds\cr
&=(e^{\lam t}-1) \int_0^\infty e^{-\lam s} D_{ij} P_s f\, ds- e^{\lam t} \int_0^t e^{-\lam s}
D_{ij} P_s f\, ds.\cr}$$ 
Since $\lam t\leq 1$, then $e^{\lam t}-1\leq c_6(\lam t)^{\al/2}$
and the bound for the second term in (3.16) now follows by using (3.15)
to bound the above integrals, and recalling again that $\lambda t\le 1$.

The proofs of (3.11) and (3.13) are similar but simpler, and are left to the 
reader (or refer to [ABP]).
\qed

\subsec{4. Existence}

Before discussing existence, we first need the following
tightness result.

\proclaim Lemma 4.1. 
Suppose $Y$ is a real-valued solution of 
$$Y_t=y_0+ M_t -\lam\int_0^t Y_r\, dr, \eqno (4.2)$$
where Ê$M_t$ is a martingale such that for some $c_1$, 
$$\angel{M}_t-\angel{M}_s\leq c_1(t-s), \qq s\leq t. \eqno (4.3)$$
Let
$T>0$, $\eps\in (0,1)$.
Let $Z_t=\int_0^t e^{-\lam (t-s)} dM_s$.
Then $Z_t=Y_t-e^{-\lam t} y_0$ and for each $q>\eps^{-1}$,
there exists a Êconstant $c_2=c_2(\eps,q,T)$ 
such that for all $\delta\in(0,1]$,
$$\E\Bigl[ \sup_{s,t\leq T, |t-s|\le \delta} |Z_t-Z_s|^{2q}\Bigr] \leq
c_2(\eps,q,T)
\frac{\delta^{\eps q-1}} {\lam^{(1-\eps)q}}. \eqno (4.4)$$

\proof 
Some elementary stochastic calculus shows that 
$$Y_t=e^{-\lam t}y_0 + \int_0^t e^{-\lam (t-s)} dM_s,$$
which proves the first assertion about $Z$.

Fix $s_0<t_0\le T$. Let
$$K_t=\Big[e^{-\lam(t_0-s_0)}-1\Big] e^{-\lam s_0}
\int_0^t e^{\lam r} dM_r$$
and
$$L_t=e^{-\lam t_0}\int_{s_0}^t e^{\lam r} dM_r.$$
Note $$Z_{t_0}-Z_{s_0}=K_{s_0}+L_{t_0}.$$
Then
$$\eqalign{\angel{K}_{s_0}
&= \Big[e^{-\lam(t_0-s_0)}-1\Big]^2e^{-2\lam s_0}
\int_0^{s_0} e^{2\lam r} d\angel{M}_r\cr
&\leq c_3 \Big[e^{-\lam (t_0-s_0)}-1\Big]^2e^{-2\lam s_0}
\frac{e^{2\lam s_0}-1}{2\lam }\cr
&\leq c_{3} Ê\Big[e^{-\lam(t_0-s_0)}-1\Big]^2 \lam^{-1}\cr
&\leq c_{3} Ê\frac{(1\land \lam(t_0-s_0))}{\lam}.\cr}$$
Considering the cases $\lam(t_0-s_0)>1$ and $\leq 1$ separately,
we see that for any $\eps\in(0,1)$ this is less than
$$c_{4}(\eps) \frac{(t_0-s_0)^\eps}{\lam^{1-\eps}}.$$
Now applying Êthe Burkholder-Davis-Gundy inequalities, we see that
$$\E |K_{s_0}|^{2q}\leq c_{5}(\eps,q) \frac{(t_0-s_0)^{\eps q}}
{\lam^{(1-\eps)q}}, Ê\qq q>1.
\eqno (4.5)$$
Similarly,
$$\eqalign{\angel{L}_{t_0}
&\leq c_6 \frac{1-e^{-2\lam(t_0-s_0)}}{2\lam}\cr
&\leq c_6\Big( \lam^{-1} \land (t_0-s_0)\Big)\cr
&= c_{6} Ê\frac{(1\land \lam(t_0-s_0))}{\lam}.\cr}$$
This leads to
$$\E |L_{t_0}|^{2q}\leq c_{7}(\eps,q) \frac{(t_0-s_0)^{\eps q}}
{\lam^{(1-\eps)q}}, Ê\qq q>1.
\eqno (4.6)$$
Combining (4.5) and (4.6) 
shows that
$$\E |Z_{t_0}-Z_{s_0}|^{2q} \leq c_{8}(\eps,q) \frac{|t_0-s_0|^{\eps q}}
{\lam^{(1-\eps) q}}.$$
It is standard to obtain (4.4) from this;
cf. the proof of Theorem I.3.11 in [B1].
\qed

Recall the definition of $H_{\beta}$ from (2.9). Ê

\proclaim Theorem 4.2. Assume $a_{ij}:H\to \R$ is continuous for all $i,j$, 
$b_i$ is continuous for all $i$, (2.2) holds, and
for some $p>1$ and positive constant $c_1$
$$\lam_k\geq c_1 k^p, \qq k\geq 1.\eqno(4.7)$$
Then for every $x_0\in H$, there is a solution $\P$ to the martingale problem for
$\sL$ starting at $x_0$. Ê
Moreover if $\beta\in(0,1)$, then any such solution has
$\sup_{\eps\leq t\leq \eps^{-1}} |X_t|_{\beta}<\infty$
for all $\eps$ Ê$\P$-a.s.
If in addition $x_0\in H_{\beta}$ for some
$\beta\in(0,1)$, then any solution $\P$ to the martingale problem for
$\sL$ starting at $x_0$ will satisfy
$$\sup_{t\le T}|X_t|_{\beta}<\infty\hbox { for all }T>0,\ \P-a.s.\eqno (4.8)$$

\proof This argument is standard and follows by making some minor
modifications to the existence result in Section 5.2 of [KX]. æWe give 
a sketch and leave the details to the reader. Fix $x_0$ in $H$.
Using the finite dimensional existence result, we may construct a solution
$X^n_t=(X^{n,k}_t:k\in\bN)$ of 
$$X_t^{n,k}=x_0(k)+1_{(k\le n)}
\Bigl[-\int_0^t\lam_kX_s^{n,k}b_k(X_s^n)ds+\sum_{j=1}^n
\int_0^t\sigma_{k,j}^n(X_s^n)dW^j_s\Bigr].$$
Here $\{W^j\}$ is a
sequence of independent one-dimensional standard Brownian motions and
$\sigma^n(x)$ is a symmetric positive definite square root of
$(a_{ij}(x))_{i,j\le n}$ which is continuous in $x\in H$ (see Lemma 5.2.1 of
[SV]). Then $X^n_t=\sum_{k=1}^n X_t^{n,k}\be_k$ has paths in
$C([0,\infty),H)$ and we next verify this sequence of processes is
relatively compact in this space. Once one has Êrelative compactness,
it is routine to use the continuity of the $a_{ij}$ and $b_i$ Êon $H$
to show that any weak limit point of $\{X^n\}$ will be a solution to the
martingale problem for $\sL$ starting at $x_0$. 

By our assumptions on $b_k$, Êeach $b_k$ is bounded above by $\gamma^{-1}$ and below
by $\gamma$. We perform a time change on $X_t^{n,k}$: let
$A^{n,k}_t=\int_0^t b_k(X_s^{n}) ds$, let $\tau^{n,k}_t$ be
the inverse of $A^{n,k}_t$, and let
$Y^{n,k}_t=X^{n,k}_{\tau^{n,k}_t}$. Then
$Y^{n,k}_t$ solves the stochastic differential equation
$$Y_t^{n,k}=x_0(k)+ 1_{(k\leq n)} \Bigl[-\int_0^t \lam_k Y^{n,k}_s ds
+M^{n,k}_t\Bigr],$$
where $M^{n,k}_t$ is a martingale satisfying Ê$|\angel{M^{n,k}}_t-
\angel{M^{n,k}}_s|\leq c_2|t-s|$, where $c_2$ is a constant not depending
on 
$n$ or $k$.

We may use stochastic calculus to write
$$Y^{n,k}_t= x^{n,k}(t) 
+Z^{n,k}_t,$$
where
$$x^{n,k}(t)=
[1_{(k\leq n)} e^{-\lam_k t} +1_{(k>n)}] x_0(k)$$
and
$$Z^{n,k}_t=1(k\le n)\int_0^t e^{-\lam_k(t-s)} dM^{n,k}_s.$$

Let $T>0$ and $s\leq t\leq T$.
Choose $\eps\in (0,1-\tfrac1{p})$ and $q>2/\eps$. Ê
By Lemma 4.1 we have for $k\le n$ and any $\delta\in(0,\gamma]$,
$$\E \Big[\sup_{u,v\leq \gamma^{-1} T, |u-v|\le
\delta\gamma^{-1}}|Z_v^{n,k}-Z_u^{n,k}|^{2q}\Big]\leq c_2(\eps,q,\gamma^{-1}T)
\gamma^{-\eps q+1}\frac{\delta^{\eps q-1}}{\lam_k^{(1-\eps)q}}.$$
Hence, undoing the time change tells us that
$$\E\Bigl[\sup_{s,t\le T, |s-t|\le \delta} |\wt X^{n,k}_t-\wt
X^{n,k}_s|^{2q}\Bigr]\leq 1_{(k\le n)}c_3(\eps,q,\gamma,T)
\frac{\delta^{\eps q-1}}{\lam_k^{(1-\eps)q}},$$
where
$$\wt X^{n,k}_t=1_{(k\le n)}(X^{n,k}_t-e^{-\lam_k \int_0^t b_k(X^n_r) dr}
x_0(k))+1_{(k>n)}x_0(k),$$
so that $\tilde X^{n,k}_{\tau_t^{n,k}}=Z_t^{n,k}$.
Now for $0\le s,t\le T$ and $|t-s|\le \gamma$,
$$\eqalign{\Big(\E |\wt X^n_t-\wt X^n_s|^{2q}\Big)^{1/q}&=
\norm{\, |\wt X^n_t-\wt X^n_s|^2}_q
=\norm{\sum_k |\wt X^{n,k}_t-\wt X^{n,k}_s|^2}_q\cr
&\leq \sum_k \norm{\, |\wt X_t^{n,k}-\wt X_s^{n,k}|^2}_q 
=\sum_k \Big(\E |\wt X^{n,k}_t-\wt X^{n,k}_s|^{2q}\Big)^{1/q}\cr
&\leq c_3(\eps,q,\gamma,T)^{1/q} \sum_k
\frac{|t-s|^{\eps-1/q}}{\lam_k^{1-\eps}},\cr}$$ 
where $\norm{\cdot}_q$ is the usual $L^q(\P)$ norm.

By our choice of $\eps$
this is bounded by $c_4(\eps,q,\gamma,T)|t-s|^{\eps/2}$,
and hence
$$\sup_n\E|\wt X^n_t-\wt X^n_s|^{2q}\leq c_4^q |t-s|^{\eps q/2}\qq
s,t\le T,\ |s-t|\le \gamma.$$ It is well known (\cite{Bi}) that this
implies the relative compactness of
$\wt X^n$ in $C(\R_+,H)$.

We may write $$X_t^n=\wt X^n_t-U^n(t),\eqno(4.9)$$ where
$$U^n(t)=\sum_{k=1}^n e^{-\lambda_k\int_0^t b_k(X_r^n)\,
dr}x_0(k)\eps_k.$$
If $s<t$, then 
$$\eqalignno{|U^n(t)-U^n(s)|^2
&=\sum_{k=1}^n\Bigl[e^{-\lambda_k\int_0^tb_k(X^n_r)\,dr}
-e^{-\lambda_k\int_0^sb_k(X^n_r)\,dr}\Bigr]^2x_0(k)^2\cr
&\le \sum_{k=1}^n ((\lambda_k^2\gamma^{-2}|t-s|^2)\wedge 1)x_0(k)^2\cr
&\le \sum_{k=1}^\infty 1_{(\lambda_k\le
\gamma|t-s|^{-1})}\lambda_k^2x_0(k)^2\gamma^{-2}|t-s|^2&(4.10)\cr
&\qquad +\sum_{k=1}^\infty 1_{(\lambda_k>\gamma|t-s|^{-1})}x_0(k)^2.\cr}
$$
Fix $\eps>0$. ÊFirst choose $N$ so that $\sum_{k=N}^\infty
x_0(k)^2<\eps$, and then $\delta>0$ so that 
$$\sum_{k=1}^\infty
1_{(\lambda_k>\gamma\delta^{-1})}x_0(k)^2<\eps,$$
and
$$\sum_{k=1}^N\lambda_k^2x_0(k)^2\gamma^{-2}\delta^2<\eps.$$
If $0<t-s<\delta$, then use the above bounds in (4.10) to conclude that 
$$\eqalign{|U^n(t)-U^n(s)|^2&\le\sum_{k=1}^N
\lambda_k^2x_0(k)^2\gamma^{-2}\delta^2+\sum_{k=N}^\infty x_0(k)^2\cr
&\qquad +\sum_{k=1}^\infty 1(\lambda_k>\gamma\delta^{-1})x_0(k)^2\cr
&<3\eps.\cr}
$$
This and the fact that $U^n(0)\to x_0$ in $H$ prove that $\{U^n\}$ is
relatively compact in $C(\R_+,H)$. ÊThe relative compactness of $\{X^n\}$
now follows from (4.9). 

Assume now $\P$ is any solution to the martingale
problem for $\sL$ starting at $x_0\in H$ and let $X^i_t$ denote
$\langle X_t,\be_i\rangle$. 
Fix $\beta\in (0,1)$ and $T>1$.
Choose $\eps\in (0, 1-\beta)$. ÊUsing a time change argument as above but
now with no parameter $n$ and $\delta=1$, we may deduce for any
$q>1/\eps$ and $k\in\bN$,
$$\eqalign{\P\Bigl(\sup_{t\le
T}|X^k_t&-e^{-\lambda_k\int_0^tb_k(X_s)\,ds}x_0(k)|
>\lambda_k^{-\beta/2}\Bigr)\cr
&\le c_5(\eps,q,T/\gamma)\lambda_k^{\beta q-q(1-\eps)}.\cr}
$$
The right-hand side is summable over $k$ by our choice of $\eps$ and
(4.7). ÊThe Borel-Cantelli lemma therefore implies that 
$$\sup_{t\le T}|X^k_t-e^{-\lambda_k\int_0^tb_k(X_s)\,ds}x_0(k)|\le
\lambda_k^{-\beta/2}\quad\hbox{ for $k$ large enough, a.s.}\eqno(4.11)$$
If $x_0\in H_\beta$, this implies that with probability $1$, for large
enough $k$,
$$\sup_{t\le T} |X_t^k|\lambda_k^{\beta/2}\le
1+x_0(k)\lambda_k^{\beta/2}\le 1+|x_0|_\beta,$$
and hence
$$\sup_{t\le T}|X_t|_\beta<\infty,\qquad a.s.$$
For general $x_0\in H$, (4.11) implies
$$\sup_{T^{-1}\le t\le T}|X_t^k|\lambda_k^{\beta/2}\le
1+e^{-\lambda_k\gamma T^{-1}}\lambda_k^{\beta/2}|x_0|\le
c_6(\gamma,T,\beta,x_0)\hbox{ for large enough $k$, Êa.s.}$$
This implies $\sup_{T^{-1}\le t\le T}|X_t|_\beta<\infty\quad a.s.$ and so
completes the proof.
\qed

\subsec{5. Uniqueness}

We continue to assume that $(a_{ij})$ and $(b_i)$ are as in Section 2 and
in particular will satisfy (2.2). ÊLet $y_0\in H$ and let $\P$ be
any solution to the martingale problem for
$\sL$ started at $y_0$. For any bounded function $f$ define
$$S_\lam f=\E \int_0^\infty e^{-\lam s} f(X_s) ds.$$
Fix $z_0\in H$ and define 
$$\sL_0 f(x)=\tfrac12 \sum_{i,j=1}^\infty a_{ij}(z_0) D_{ij} f(x)
-\sum_{i}^\infty
\lam_ix_i b_i(z_0)D_if(x). \eqno (5.1)$$
Set $\sB=\sL-\sL_0$ and let $R_\lam$ be the resolvent for $\sL_0$
as Êin Section 2.

ÊTo make this agree with the definition of $\sL_0$ in
Section~2 we must replace $\lambda_i$ by
$\wh\lambda_i=b_i(z_0)\lambda_i$ and set $a_{ij}^0=a_{ij}(z_0)$. ÊAs
$\gamma\le b_i(z_0)\le \gamma^{-1}$, and the constants in Corollary 3.5
may depend on $\gamma$, we see that the bounds in Corollary 3.5
involving the original $\lambda_i$ remain valid for $R_\lambda$. ÊWe
also will use the other results in Section 3 with $\wh\lambda_i$ in
place of $\lambda_i$ without further comment. ÊIn addition, if we
simultaneously replace $b_i$ by $\wh b_i=b_i/b_i(z_0)$, then 
$$\sL f(x)={1\over 2}\sum_{i,j=1}^\infty
a_{ij}(x)D_{ij}f(x)-\sum_{i=1}^\infty\wh\lambda_ix_i\wh b_i(x)D_if(x),$$
$$\sL_0f(x)={1\over 2}\sum_{i,j=1}^\infty
a_{ij}(z_0)D_{ij}f(x)-\sum_{i=1}^\infty\wh\lambda_ix_iD_if(x),$$
and
$$\wh b_i(z_0)=1\quad\hbox{all }i.$$
In Propositions 5.1 and 5.2 below we will simply assume $b_i(z_0)=1$
for all $i$ without loss of generality, Êit being understood that the above
substitutions are being made. ÊIn each case it is easy to check that the
hypotheses on $(b_i,\lambda_i)$ carry over to $(\wh b_i,\wh\lambda_i)$
and as the conclusions only involve $\sL$, $\sL_0$, $R_\lambda$, and our
solution $X$, which remain unaltered by these substitutions, this
reduction is valid.

ÊLet
$$\eta=\sup_x\sum_{i,j=1}^\infty |a_{ij}(x)-a_{ij}(z_0)|.\eqno(5.2)$$
Set $$B_i(x) =x_i(b_i(x)-1).$$ ÊAs before, $\alpha$ will denote a
parameter in $(0,1)$.

\proclaim Proposition 5.1. Assume
$$\sum_{i\leq j} |a_{ij}|_{C^\al} \lam_j^{-\al/2}<\infty, \eqno (5.3)$$
$$\sum_i \lam_i^{1/2} \norm{B_i}_{C_b}<\infty, \eqno (5.4)$$
and
$$\sum_i \lam_i^{(1-\al)/2} |B_i|_{C^\al}<\infty. \eqno (5.5)$$
There exists $c_1(\lambda)\to 0$ as $\lam\to \infty$ and
$c_2=c_2(\al,\gamma)$ such that
for all $f\in C^\al$, we have $\sB R_\lam f\in C^\al$ and
$$\norm{\sB R_\lam f}_{C^\al}\leq (c_1(\lam)+c_2\eta)\norm{f}_{C^\al}.$$

\proof 
We have
$$\eqalignno{|\sB R_\lam f(x)|&\leq \sum_{i,j} |a_{ij}(x)-a_{ij}(z_0)|\, |D_{ij}
R_\lam f(x)|\cr
&\qq\qq + \sum_i \lam_i|x_i|\, |b_i(x)-1|\, |D_iR_\lam f(x)|\cr
&\leq \eta c_3 Ê|f|_{C^\al} 
+ c_4(\lambda)|f|_{C^\al},&(5.6)\cr}$$
where $c_4(\lambda)\to 0$ as $\lambda\to\infty$ by (5.4) and (3.11). ÊIn
particular, the series defining $\sB R_\lam f$ is absolutely
uniformly convergent.

Let $\wh a_{ij}(x)=a_{ij}(x)-a_{ij}(z_0)$. If $h\in H$, then
$$\eqalignno{|\sB ÊR_\lam f(x+h)-\sB R_\lam f(x)|&
=\Bigl|\sum_{i,j} [\wh a_{ij}(x+h)D_{ij}R_\lam f(x+h)-\wh a_{ij}(x) D_{ij} R_\lam
f(x)]\cr
&\qq\qq + \sum_i \lam_i[B_i(x+h)D_iR_\lam f(x+h)-B_i(x)D_iR_\lam Êf(x)]\Bigr|\cr
&\leq \Bigl| \sum_{i,j} \wh a_{ij}(x+h) (D_{ij} R_\lam f(x+h)-D_{ij} R_\lam f(x)
)\Bigr|\cr
&\qq +\Bigl|\sum_{i,j} (\wh a_{ij}(x+h)-\wh a_{ij}(x))D_{ij} R_\lam f(x) \Bigr|\cr
&\qq\qq +\Bigl|\sum_i \lam_i B_i(x+h) (D_iR_\lam f(x+h)-D_iR_\lam f(x))\Bigr|\cr
&\qq\qq+ \Bigl|\sum_i \lam_i (B_i(x+h)-B_i(x)) D_iR_\lam f(x)\Bigr|\cr
&=S_1+S_2+S_3+S_4.&(5.7) \cr}$$
Use (3.14) to see that
$$\eqalignno{S_1&\leq c_5\sum_{i,j} |\wh a_{ij}(x+h)| Ê|f|_{C^\al} |h|^\al\cr
&\leq c_6\eta Ê|f|_{C^\al} |h|^\al.&(5.8)\cr}$$
By (3.12)
$$\eqalignno{S_2&\leq \sum_{i,j} |a_{ij}(x+h)-a_{ij}(x)|\, |D_{ij} R_\lam f(x)|\cr
&\leq c_7\sum_{i\leq j} |a_{ij}|_{C^\al} |h|^\al Ê(\lam+\lam_j)^{-\al/2}
Ê|f|_{C^\al}\cr
&\leq c_{8}(\lam) |f|_{C^\al} |h|^\al, &(5.9)\cr}$$
where (5.3) and dominated convergence imply $\lim_{\lam\to \infty} c_{8}(\lam)=0$.
By (3.13)
$$S_3\leq c_9\sum_i \lam_i |B_i(x+h)| Ê(\lam+\lam_i)^{-1/2} |f|_{C^\al} |h|^\al
\leq c_{10}(\lam) |f|_{C^\al} |h|^\al, \eqno (5.10)$$
where $c_{10}(\lam)\to 0$ as $\lam\to \infty$ by (5.4) and dominated
convergence. By (3.11)
$$S_4\leq c_{11}\sum_i \lam_i |B_i|_{C^\al} Ê(\lam+\lam_i)^{-(1+\al)/2} |f|_{C^\al}
|h|^\al\leq c_{12}(\lam) |f|_{C^\al} |h|^\al, \eqno (5.11)$$
where again $c_{12}(\lam)\to 0$ as $\lam\to \infty$ by (5.5). 
Combining (5.8), (5.9), (5.10), and (5.11) yields
$$|\sB R_\lam f|_{C^\al} \leq [c_{13}(\lam)+c_{14} \eta]|f|_{C^\al}.$$
This and (5.6) complete the proof.
\qed

Let $C^\al_n$ denote those functions in $C^\al$ which only depend on the first
$n$ coordinates. Note that $\sT_0\subset\cup_n C^\al_n$.
\proclaim Proposition 5.2. If $f\in \cup_n C^\al_n$, then
$$S_\lam f=R_\lam f(y_0)+S_\lam \sB R_\lam f. \eqno (5.12)$$

\proof
Fix
$z_0\in H$. Suppose $h\in \sT$. Since $h(X_t)-h(X_0)-\int_0^t \sL h(X_s)
ds$ is a martingale, taking expectations we have
$$\E h(X_t)-h(y_0)=\E\int_0^t \sL h(X_s) ds.$$
Multiplying by $e^{-\lam t}$ and integrating over $t$ from
0 to $\infty$, we obtain
$$S_\lam h -\frac{1}{\lam } h(y_0)=\E\int_0^\infty e^{-\lam t} \int_0^t
\sL h(X_s) ds \, dt
=\frac{1}{\lam} \E \int_0^\infty e^{-\lam s} \sL h(X_s) ds
=\frac{1}{\lam} S_\lam \sL h.$$
This can be rewritten as
$$\lam S_\lam h- S_\lam \sL_0 h= h(y_0)+S_\lam \sB h. \eqno (5.13)$$

Define 
$$\sL_0^n f(x)=\sum_{i,j=1}^n a_{ij}(z_0) D_{ij}f(x)-\sum_{i=1}^n
\lam_ix_i D_if(x).$$
Let $R_\lam^n$ be the corresponding resolvent. The corresponding process
is an $n$-dimensional Ornstein-Uhlenbeck process which starting from $x$
at time $t$ is Gaussian with mean vector $(x_ie^{-\lam_it})_{i\le n}$ and
covariance matrix
$C_{ij}(t)=a_{ij}(z_0)(1-e^{-(\lam_i+\lam_j)t})(\lam_i+\lam_j)^{-1}$.
These parameters are independent of $n$ and the distribution coincides
with the law of the first $n$ coordinates (with respect to ${\be_i}$) of
the process with resolvent $R_\lam$. 

Now take $f\in C^\al_n$ and let $h(x)=R_\lam
f(x)=R_\lam^nf(x_1,\dots,x_n)$. (Here we abuse our notation slightly by
having $f$ also denote its dependence on the first $n$ variables.) By
Corollary 3.5 and (3.10),
$h\in\sT$. Moreover, 
$\sL_0h=\sL_0^nR^n_\lam f=\lam R^n_\lam f-f=\lam R_\lam f-f$.
The second equality is standard since on functions in $C_b^2$, $\sL_0^n$
coincides with the generator of the finite-dimensional diffusion. Now
substitute this into (5.13) to derive (5.12). 
\qed

To iterate (5.12) we will need to extend it to $f\in C^\al$ by an
approximation argument. ÊRecall $\wh\lambda_i=b_i(z_0)\lambda_i$. 
\ms
\ni{\bf Notation.} Write $f_n\bpto f$ if $\{f_n\}$ converges to $f$
pointwise and boundedly.

\proclaim Lemma 5.3. (a) 
ÊÊIf $f\in C^\al$, then $pR_pf \bpto f$
as $p\to \infty$ and
$$\sup_{p>0} \norm{pR_p f}_{C^\al} \leq \norm{f}_{C^\al}.$$
\item{(b)} For $p>0$ there is a $c_1(p)$ such that for any bounded
measurable
$f:H\to\R$, $R_pf\in C^\al$ and
$\norm{pR_pf}_{C^\al}\leq c_1(p) \norm{f}_{C_b}$.

\proof (a)
Note if $f\in C^\alpha$, then
$$\norm{pR_p f}_{C_b}\leq \int_0^\infty pe^{-pt} \norm{P_tf}_{C_b} dt
\leq \norm{f}_{C_b}$$
and
$$pR_pf(x)-f(x)=\int_0^\infty pe^{-pt} (P_tf(x)-f(x)) dt\to 0$$
because $P_tf(x)\bpto f(x) $ as $t\to 0$.

Let $X_t$ be the solution to (2.4) (so that $X$ has resolvents
$(R_\lambda)$) and let\hfil\break
$X_t^i=\angel{X_t, \be_i}\be_i$. Then $X_t^i$
satisfies
$$X_t^i=X_0^i+M_t^i -\wh\lam_i \int_0^t X_s^i ds, \eqno (5.14)$$
where $M_t^i$ is a one-dimensional Brownian motion with $\Cov(M_t^i, M_s^i)=
a_{ii}(s\land t)$. Let $X_t^{x_i,i}$ denote the solution to (5.14) when
$X_0^i=x_i$. Then
$$X_t^{x_i+h_i,i}-X_t^{x_i,i}=h_i-\wh\lam_i \int_0^t
(X_s^{x_i+h_i,i}-X_s^{x_i,i})ds,$$ and so
$$X^{x_i,h,i}_t-X^{x,i}_t=e^{-\wh\lam_i t} h_i \be_i.$$
Hence, if $X_t^x$ is defined by $\angel{X_t^x, \be_i}=X_t^{x_i,i}$,
$$|X^{x+h}_t-X^x_t|=\Bigl|\sum h_i^2 e^{-2\wh\lam_i t}\Bigr|^{1/2}\leq
|h|.$$ Therefore
$$|P_tf(x+h)-P_tf(x)|\leq |f|_{C^\al} \E(|X^{x+h}_t-X^x_t|^\al)\leq |f|_{C^\al} |h|^\al$$
and so
$$|pR_pf(x+h)-pR_pf(x)|\leq \int_0^\infty pe^{-pt} |P_tf(x+h)-P_tf(x)|dt\leq
|f|_{C^\al} |h|^\al,$$
i.e., $|pR_pf|_{C^\al}\leq |f|_{C^\al}$. This proves (a).

(b) As we mentioned above, for any bounded measurable $f$, æ$\norm{pR_p
f}_{C_b}\leq
\norm{f}_{C_b}$. We also have Ê Ê
$$\eqalign{P_s pR_pf-pR_pf&=\int_0^\infty pe^{-pt} [P_{s+t}f-P_tf] dt\cr
&=(e^{ps}-1) \int_0^\infty pe^{-pt} P_t f\, dt-e^{ps}\int_0^s
pe^{-pt} P_t f\, dt.\cr}$$
The right hand side is bounded by
$$2(e^{ps}-1) \norm{f}_{C_b}.$$
This in turn Êis bounded by $c_2(p) s^{\al/2}$ for $0\le s\le 1$.
Also,
$$\norm{P_spR_pf-pR_pf}_{C_b}\leq 2\norm{f}_{C_b}\leq 2s^{\al/2} 
\norm{f}_{C_b}\quad \hbox{ for $s\geq 1$}.$$
ÊHence $\norm{pR_p f}_\bal\leq c_3(p) \norm{f}_{C_b}$.
Our conclusion follows by (3.10), which holds for the $\{\wh\lambda_i\}$
just as it did for $\{\lambda_i\}$.
\qed

\proclaim Lemma 5.4. Suppose $f_n\bpto 0$ where
$\sup_n\norm{f_n}_{C^\al}<\infty$. Then 
$$D_{ij} R_\lam f_n\bpto 0\hbox{ and }
D_iR_\lam f_n\bpto 0\hbox{ as }n\to\infty\hbox{ for all }i,j.$$

\proof We focus on the second order derivatives as the proof for the first
order derivatives is simpler. ÊWe know from Corollary 3.3 that
$D_{ij} R_\lam f_n$ is uniformly bounded in
$C^\al$ norm, so in particular, it is uniformly bounded in ${C_b}$ norm and
we need only establish the pointwise convergence. We have from (3.8) that 
$$\norm{D_{ij} P_t f_n}_{C_b}\leq c_1 \norm{f_n}_{C^\al} t^{\al/2
-1}.\eqno(5.15)$$ 
{}From Proposition 3.2, we have
$$D_{ij}P_tf_n=D_iP_{t/2}D_{Q_{t/2}\be_j}P_{t/2}f_n.\eqno(5.16)$$
Fix $t>0$ and $w\in H$. The proof of Proposition 5.2 in [ABP] shows there exist 
random variables $R(t,w)$ and $Y_t$ such that
$$D_wP_tf(x)=\E[f(Q_tx+Y_t)R(t,w)], \qq f\in {C_b}(H),$$
and
$$\E [R(t,w)^2]\leq \frac{|w|^2}{\gamma t}.$$
Therefore
$$
h_n(j,t,x)\equiv D_{Q_{t/2}e_j}P_{t/2}f_n(x)=\E(f_n(Q_{t/2}x+Y_{t/2})R(t/2,Q_{t/2}\be_j))
\bpto 0$$
by Êdominated convergence. Moreover Ê
Cauchy-Schwarz implies
$$\norm {h_n(j,t)}_{C_b}\le Ê(\gamma t)^{-1/2} \sup_m\norm{f_m}_{C_b}.$$
ÊRepeating the above reasoning and using (5.16) we have
$$D_{ij}P_tf_n(x)=D_iP_{t/2}h_n(x)=\E(h_n(Q_{t/2}x+Y_{t/2})R(t/2,\be_i))\bpto 0
,$$
and
$$\norm{D_{ij}P_tf_n}_{C_b}\le (\gamma
t)^{-1}\sup_m \norm{f_m}_{C_b}.\eqno(5.17)$$
Fix $\eps>0$. Write
$$|D_{ij}R_\lam f_n(x)|\le\Bigl|\int_0^\eps e^{-\lam t} D_{ij}P_t f_n(x)
dt\Bigr| +\Bigl|\int_\eps^\infty e^{-\lam t} D_{ij}P_t f_n(x) dt\Bigr|;$$
by dominated convergence and (5.17) the second term tends to 0,
while (5.15) shows the first term is bounded by
$$\int_0^\eps c_2 \norm{f_n}_{C^\al} t^{\al/2 -1}dt\leq c_3(\sup_m
\norm{f_m}_{C^\al}) \eps^{\al/2}.$$ Therefore
$$\limsup_{n\to \infty} |D_{ij}R_\lam f_n(x)| \leq c_4(\sup_m
\norm{f_m}_{C^\al}) \eps^{\al/2}.$$ Since $\eps$ is arbitrary, 
$$\limsup_{n\to \infty} |D_{ij}R_\lam f_n(x)|=0.$$
\qed

\proclaim Proposition 5.5. Assume (5.4). ÊIf $f\in C^\al$, then
$$S_\lam f=R_\lam f(y_0)+ S_\lam \sB R_\lam f.\eqno (5.18)$$

\proof We know $f_p=f-pR_pf\bpto 0$ as $p\to\infty$ by Lemma 5.3. Lemma
5.3 also shows
$\norm{f_p}_{C^\al}\le 2\norm{f}_{C^\al}$, and therefore we may use Lemma 5.4, the
finiteness of $\eta$, (5.4) (in fact a weaker condition suffices here),
and dominated convergence to conclude
$$\eqalign{\sB R_\lam f_p(x)=\sum_{i,j}&(a_{ij}(x)-a_{ij}(z_0))D_{ij}(R_\lam f_p)(x)\cr
& +\sum_i \lam_i x_i(b_i(x)-b_i(z_0)) D_i(R_\lam f_p)(x)\bpto 
0\hbox{ as $p\to\infty$}.\cr}$$
Here we also use the bounds $\norm{D_{ij}R_\lam f_p}_{C_b}\le c\Vert
f\Vert_{C^\al}$ and
$\norm{D_iR_\lam f_p}_{C_b}\le c\lam_i^{-1/2}\norm{f}_{C^\al}$ from (3.11),
(3.12) and Lemma 5.3(a). By using dominated convergence it is now
easy to take limits through the resolvents to see that to prove (5.18) it
suffices to fix $p>0$ and verify it for $f=pR_ph$ where $h\in C^\al$. Fix
such an $h$.

Let $z_n(x)=\sum_{i=1}^nx_i\be_i+\sum_{i>n}(z_0)_i\be_i\to x\hbox{ as
}n\to\infty$ and define $h_n(x)=h(z_n(x))$. Then $h_n\bpto h$ since $h\in
C^\al$. Recall the definition of $R^n_p$ from the 
proof of Proposition 5.2;
by the argument there, we see that the function 
$pR_ph_n(x)=pR_p^nh_n(x_1,\dots,x_n)$ depends only on
$(x_1,\dots,x_n)$. By Lemma~5.3(b) $pR_ph_n\in C^\al$ and therefore is in
$C_n^\al$. Proposition 5.2 shows that (5.18) is valid with $f=R_ph_n$.
Now $pR_ph_n\bpto pR_ph$ as $n\to\infty$ and $\sup_n\norm{pR_ph_n}_{C^\al}\le
c_1(p)$ by Lemma 5.3(b). Therefore if $d_n=pR_p(h_n-h)$ we may use Lemma
5.4, Corollary 3.5, and dominated convergence, as before, to conclude
$$\eqalign{\sB R_\lam
d_n(x)=&\sum_{i,j}(a_{ij}(x)-a_{ij}(z_0))D_{ij}(R_\lam d_n)(x)
\cr
\quad&+\sum_i\lam_ix_i(b_i(x)-b_i(z_0))D_i(R_\lam d_n)(x)
\bpto 0
\hbox{ as }n\to\infty.\cr}$$ 
We may now let $n\to\infty$ in (5.18) with
$f=pR_ph_n$ to derive (5.18) with $f=pR_ph$, as required.\qed

\proclaim Theorem 5.6. Assume (2.2), Êeach $a_{ij}$ and each $b_i$ is
continuous, Ê(4.7), Ê(5.3), (5.4), and (5.5) Êhold.
There exists $\eta_0$, depending only on
$(\al,\gamma)$, such that if 
$\eta\leq
\eta_0$, then for any $y_0\in H$ there is a unique solution to the
martingale problem for $\sL$ started at $y_0$.

\proof Existence follows from Theorem 4.2.

Let $\P$ be any solution to the martingale
problem and define $S_\lam $ as above. Suppose $f\in C^\al$. Then by
Proposition 5.5 we have
$$S_\lam f=R_\lam f(y_0)+ S_\lam \sB R_\lam f.$$
Using Proposition 5.1 we can iterate the above and obtain
$$S_\lam f=R_\lam\Big(\sum_{i=0}^k (\sB R_\lam )^i\Big)f(y_0)+ S_\lam (\sB R_\lam)^{k+1} f.$$
Provided $\eta_0=\eta_0(\al,\gamma)$ is small enough, our hypothesis that
$\eta\leq
\eta_0$ and Lemma 5.1 imply
that for $\lam>\lam_0(\al,\gamma,(a_{ij}), (b_i))$, 
the operator $\sB R_\lam $ is bounded on $C^\al$ with norm
strictly less than $\frac12$. Therefore $\sum_{i=k+1}^\infty (\sB R_\lam
)^i f$ converges to 0 and $(\sB R_\lam )^{k+1} f$
also converges to 0, both in $C^\al$ norm, as $k\to \infty$. In
particular, they converge to 0 in sup norm, so
$R_\lam(\sum_{i=k+1}^\infty (\sB R_\lam)^i )f(y_0)$
and $S_\lam (\sB R_\lam)^{k+1} f$ both converge to 0 as $k\to\infty$.
It follows that
$$S_\lam f=R_\lam \Big(\sum_{i=0}^\infty(\sB R_\lam)^i \Big)f(y_0).
$$
This is true for any solution to the martingale problem, so
$S_\lam$ is uniquely defined for large enough $\lam$. Inverting
the Laplace transform and using the
continuity of $t\to \E f(X_t)$, we see that for every $f\in C^\al$, $\E f(X_t)$ has
the same value for every solution to the martingale problem. It is not
hard to see that $\sT_0\subset C^\al$ is dense with respect to the topology
of bounded pointwise convergence in the set of all bounded functions.
{}From here standard arguments (cf. [B2], Section æVI.3) allow us to
conclude the uniqueness of the martingale problem of $\sL$ starting
at $y_0$ as long as we have $\eta\leq \eta_0$.
\qed

Set
$$Q_{\beta,N}=\{x\in H:|x|_{\beta}\leq N\}.$$

\proclaim Theorem 5.7. Assume $(b_i)$ and $(a_{ij})$ are as in Section 2,
so that (2.2) holds. ÊAssume also that $\alpha,\beta\in(0,1)$ satisfy
\item{(a)} There exists $p>1$ and $c_1>0$ such that $\lam_j\geq c_1j^p$.
\item{(b)} $\sum_{i\leq j} |a_{ij}|_{C^\al} \lam_j^{-\al/2}<\infty$.
\item{(c)}$\sum_j \lam_j^{-\beta}<\infty$.
(For example, this holds if $\beta>1/p$.)
\item{(d)} For all $N>0$, for all $\eta_0>0$, and for all $x_0\in Q_{\beta,N}$ there exists $\delta>0$ such that if
$|x-x_0|<\delta$ and $x\in Q_{\beta,N}$, then
$$\sum_{i,j} |a_{ij}(x)-a_{ij}(x_0)|<\eta_0.$$
\item{(e)} $\sum_i \lam_i^{1/2} |b_i|_{C^\al}<\infty$.
\par
\ni {\sl Then for all $y\in H_{\beta}$ there exists a unique solution to the 
martingale problem for $\sL$ starting at $y$.
}

\bs
\noindent{\bf Remark.} By Theorem 4.2, any solution to the martingale
problem for
$\sL$ starting at $y\in H$ will immediately enter $H_\beta$ and remain
there a.s. for any $\beta\in(0,1)$. ÊHence the spaces $H_\beta$ are
natural state spaces for the martingale problem.

\proof
Fix $\beta\in (0,1)$ as in (c) and write $Q_N$ for $Q_{\beta,N}$. Let
$\P$ be a solution to the martingale problem for $\sL$. By Theorem 4.2
we
only need consider uniqueness. If $T_N=\inf\{t: X_t\notin Q_N\}$, 
then by Theorem 4.2 we see
that
$T_N\uparrow \infty $, a.s. and it suffices to show uniqueness for
$\P(X_{\cdot\land T_N}\in \cdot)$. (c) implies $Q_N$ is compact and
so as in the proof of Theorem VI.4.2 of [B2] it suffices to show:
\ms

\item{(5.19)} for all $x_0\in Q_N$ there exist $r>0$, $\wt a_{ij}$, and
$\wt b_i$ such that
$a_{ij}=\wt a_{ij}$ and $b_i=\wt b_i$ on $Q_N\cap \{x\in H: |x-x_0|<r\}$ and the martingale problem for $\wt \sL$ starting at $y$ has a unique solution for all $y\in Q_N$.
Here $\wt \sL$ is defined analogously to $\sL$ but with $a_{ij}$ and $b_i$ replaced
by $\wt a_{ij}$ and $\wt b_i$, respectively.
\ms

Fix $x_0\in Q_N$, $\eta_0$ as in Theorem 5.6. 
Choose $\delta$ as in (d).
We claim we can choose $1\ge\delta_1>0$ depending on $\delta$ and
$N$ Êsuch that if $x\in Q_N$ and $\norm{x-x_0}_\infty <\delta_1$, then
$|x-x_0|<\delta$. Here $|x|_\infty=\sup_i |\langle x,\be_i\rangle |$.

To prove the claim, note that $\norm{x-x_0}_\infty\leq \delta_1$
implies
that for any $K_0$
$$\sum_k (x^k -x^k_0)^2\leq \sum_k \delta_1^2\land (4N^2\lam_k^{-\beta})
\leq K_0\delta_1^2 +4N^2 \sum_{k>K_0} \lam_k^{-\beta}.$$
So first choose $K_0$ such that the second term is less than
$\delta^2/2$ and then set $\delta_1=\delta/\sqrt{2K_0}$. 

Now let $[p_j,q_j]=[x_0^j-\delta_1,x_0^j+\delta_1]\cap [-N\lam_j^{-\beta/2},
N\lam_j^{-\beta/2}]$ and note
$p_j<q_j$ as $x_0\in Q_N$. Let $\psi_j:\R\to \R$ be defined
by 
$$\psi_j(x)=\cases{x &if $p_j\leq x\leq q_j$;\cr
p_j& if $x<p_j$;\cr Êq_j& if $x>q_j$.\cr}$$
Define $\psi:H\to Q_N\cap \{x\in H: \norm{x-x_0}_\infty<\delta_1\}$
by
$$\psi(x)=\sum_{j=1}^\infty \psi_j(\angel{x,e_j})e_j.$$
As $\norm{\psi_j}_\infty^2 \leq N^2\lam_j^{-\beta}$, $\psi$ is well
defined by (c).

Take $r=\delta_1\in(0,1]$ and set $\wt a_{ij}(x)=a_{ij}(\psi(x))$. 
If $|x-x_0|<r$ and $x\in Q_N$, then Ê$\norm{x-x_0}_\infty<r$ and therefore
$\psi(x)=x$, which says that $\wt a_{ij}(x)=a_{ij}(x)$ for all $i, j$.

Define 
$$\rho(u)=\cases{u&if $|u|<r$;\cr 
(2r-|u|)u/r&if $r\leq |u|<2r$;\cr 0&if $2r\leq |u|$,\cr}$$
and set $\wt b_i(x)=b_i(x_0+\rho(x-x_0))$. ÊIf $|x-x_0|<r$, then
$\rho(x-x_0)=x-x_0$ and so $\wt b_i(x)=b_i(x)$. ÊAlso $\wt b_i$ is
clearly continuous as (e) implies that $b_i$ is.

We now show that $\wt a_{ij}$ satisfies the hypotheses of Theorem 5.6.
For any $x$
$$\sum_{i,j} |\wt a_{ij}(x)-\wt a_{ij}(x_0)|=\sum_{i,j} |a_{ij}(\psi(x))
-a_{ij}(x_0)|. \eqno (5.20)$$
Since $\norm{\psi(x)-x_0}_\infty\leq r$ and $\psi(x)\in Q_N$, it
follows that Ê
$|\psi(x)-x_0|<\delta$. (d) now implies that the right hand side of (5.20)
is less than $\eta_0$. It remains only to check (5.3)
for $\wt a_{ij}$. But
$$|\psi_j(x)-\psi_j(x+h_j)|\leq |h_j|,$$
and so 
$$|\psi(x)-\psi(x+h)|\leq |h|.$$
Therefore
$$\eqalign{|\wt a_{ij}(x+h)&-\wt a_{ij}(x)|\cr
&=|a_{ij}(\psi(x+h))-a_{ij}(\psi(x))|\cr
&\leq |a_{ij}|_{C^\al} |\psi(x+h)-\psi(x)|^\al\cr
&\leq |a_{ij}|_{C^\al} |h|^\al,\cr}$$
and so 
$$|\wt a_{ij}|_{C^\al}\leq |a_{ij}|_{C^\al}.$$
Hence $\wt a_{ij}$ satisfies (5.3) because $a_{ij}$ does.

If we set $B_i(x)=x_i(\wt b_i(x)-\wt b_i(x_0))$,
it is easy to check that $B_i(x)$ is 0 for $|x-x_0|\geq 2r$,
$\norm{B_i}_\infty\leq c_{2} |\wt b_i|_{C^\al}\le c_2|b_i|_{C^\al}$, and
$|B_i|_{C^\al}\leq c_{2} |\wt b_i|_{C^\al}\le c_2|b_i|_{C^\al}$, where $c_{1}$
may depend on
$x_0$.
Therefore (e) implies $(\wt b_i)$ satisfies (5.4) and (5.5). 

We see then that Theorem 5.6 applies to $\wt a_{ij}$ and $\wt b_i$ and so (5.19)
holds. \qed

\bs
\noindent
{\bf Example 5.8.} We discuss a class of examples where the $b_i=1$ and the $a_{ij}$ are
zero unless $i$ and $j$ are sufficiently close together. Let $M\in\bN$, $\al\in(0,1)$ and $S_M(i,j)$ be the
subspace of $H$ generated by $\{\be_k:|k-i|\vee |k-j|\le M\}$. ÊAlso let
$\Pi_{S_M(i,j)}$ be the projection operator onto $S_M(i,j)$. ÊAssume that
$a_{ji}(x)=a_{ij}(x)=\langle\be_i,a(x)e_j\rangle$ satisfies (2.2) and
depends only on coordinates corresponding to $S_M(i,j)$, that is, 
$$a_{ij}(x)=a_{ij}(\Pi_{S_M(i,j)}x)\quad\hbox{ for all }x\in H,\quad
i,j\in\bN.\eqno(5.21)$$
In particular, (5.21) implies $a_{ij}$ is constant if $|i-j|>2M$.
Also suppose that
Ê$$\sup_{i,j}|a_{ij}|_{C^\al}=c_1<\infty.\eqno(5.22)$$
ÊÊSet
$b_i(x)=1$ for all $i$, $x$ and also assume
$$
\lam_j\ge c_2 j^p\quad\hbox{all }j\hbox{ for some }p>1,\eqno(5.23)$$
and $\beta\in(0,1)$ satisfies
$$\sum_{j=1}^\infty\lam_j^{{-\beta\al\over
2}+\delta}<\infty\quad\hbox{for some }\delta>0.\eqno(5.24)
$$
For example, (5.24) will hold if $p>2$ and $\beta\alpha>2/p$. We then
claim that the hypotheses of Theorem 5.7 hold and so there is a unique
solution to the martingale problem for
$\sL f(x)=\sum_{i,j}a_{ij}(x)D_{ij}f(x)-\sum_i\lam_ix_i D_if(x)$, starting
at any $y\in H_\beta$.

We must check conditions (b)--(d) of Theorem 5.7. ÊNote first that 
$$|a_{ij}(x+h)-a_{ij}(x)|\le 1_{(|i-j|\le 2M)}|a_{ij}|_{C^\al} |h|^\al,$$
so that $|a_{ij}|_{C^\al}\le 1_{(|i-j|\le 2M)}c_3$ and hence by (5.24),
$$\sum_{i\le
j}|a_{ij}|_{C^\al}\lam_j^{-\al/2}\le(2M+1)c_5\sum_j\lam_j^{-\al/2}<\infty.$$
This proves (b), and (c) is immediate from (5.24). If $N>0$, $x,x_0\in
Q_{\beta,N}$, then for small enough $\eps>0$,
$$\eqalign{
\sum_{i,j}|a_{ij}(x)-a_{ij}(x_0)|
&\le 2\sum_{i\le
j}|a_{ij}|_{C^\al}\Bigl[\sum_k1(|k-i|\vee|k-j|\le
M)(x(k)-x_0(k))^2\Bigr]^{\al/2}\cr
&\le 2|x-x_0|^\eps \sum_i \sum_{j=1}^{i+2M} \Big[\sum_k 1_{(|k-i|\leq M)} 
|x(k)-x_0(k)|^{2-(2\eps/\al)}\Big]^{\al/2}\cr
&\le |x-x_0|^\eps c_4(M)\sum_{k=1}^\infty |x(k)-x_0(k)|^{\al-\eps}\cr
&\le c_5(M)|x-x_0|^\eps \sum_{k=1}^\infty
(2N)^{\al-\eps}\lam_k^{{-\beta\over 2}(\al-\eps)}\cr
&\le c(M,N)|x-x_0|^\eps.\cr}
$$
We have used (5.22), $x,x_0\in Q_{\beta,N}$ and (5.24) in the above. ÊThis proves
(d), as required.

\medskip

\ni {\bf Example 5.9.} We give a more specific realization of the previous example.
Continue to assume $b_i=1$ for all $i$, (5.23), and (5.24). Let
$L,N\geq 1$ (we can take $N=1$, for example)
and for $k\geq 1$ let $I_k=\{(k-1)N+1, \ldots, kN\}$.
For each $k$ assume $a^{(k)}: {\R}^{2L+N}\to \sS^+_N$, the space of
symmetric positive definite $N\times N$ matrics. Assume
for all $k$, for all $x\in \R^{2L+N}$, and for all $z\in \R^N$,
$$\sum_{i=1}^N \sum_{j=1}^N a^{(k)}_{ij}(x)z_iz_j \in [\gamma |z|^2, \gamma^{-1} |z|^2],
Ê\eqno (5.25)$$
and
$$\sup_k \max_{1\leq i,j\leq N} |a^{(k)}_{ij}|_{C^\al}<\infty. \eqno (5.26)$$
Now for $x\in H$, let $\pi_k x=(\langle x, \be_{((\ell+k-1)N-L)\lor 1}\rangle)_{\ell =1, \ldots, 2L+N}
\in \R^{2L+N}$
and define $a: H\to L(H,H)$ by
$$\eqalign{\langle a(x)\be_i, \be_j\rangle &=a_{ij}(x)=a_{ji}(x)\cr
&=\cases{a^{(k)}_{i-(k-1)N, j-(k-1)N}(\pi_kx)&if $i,j\in I_k, k\geq 1$,\cr
0&if $(i,j)\notin \cup_{k=1}^\infty I_k\times I_k$.\cr}\cr}$$
Then for all $x,z\in H$,
$$\eqalign{\sum_i\sum_j a_{ij}(x)z_iz_j&=\sum_{k=1}^\infty\sum_{i,j\in I_k} a_{ij}(x)z_iz_j\cr
&=\sum_{k=1}^\infty \sum_{i,j=1}^N a_{ij}^{(k)} (\pi_kx)z_{(k-1)N+i}z_{(k-1)N+j}\cr
&\in [\gamma|z|^2, \gamma^{-1}|z|^2] \cr}$$
by (5.25), and so (2.2) holds. Note that if $i,j\in I_k$, then (using the notation
of Example 5.8)
$\sS_{L+N}(i,j)\supset \{(k-1)N-L+1, \ldots kN+L\}$, and so
(5.21) with $M=L+N$ is immediate from the above definitions. Also (5.22) is implied by (5.26).
The conditions of Example 5.8 therefore hold and so weak existence and uniqueness of solutions
hold for the martingale problem for $\sL$ with initial conditions in $H_\beta$.

\medskip

\ni {\bf Remark 5.10.} The above examples demonstrate the novel features of our results.  The fact
that our perturbation need not be non-negative facilitates the localization argument (see Remark 9
in [Z] for comparison) and the presence of $\{\lambda_j^{-\alpha/2}\}$ in condition (b) of
Theorem~5.7 means that the perturbation need not be H\"older in the trace class norm.  The latter
allows for the possibility of locally dependent H\"older coefficients with just bounded H\"older
norms, something that seems not to be possible using other results in the literature. On the other
hand [Z] includes an SPDE example which our approach cannot handle in general unless, for
example, the orthonormal basis in the equation diagonalizes the second derivative operator.  This
is because he has decoupled the conditions on the drift operator  and noise term, while ours are
interconnected.  The latter leads to the double summation in conditions (b) and (d) of Theorem~5.7,
as opposed to the trace class conditions in [Z].  All of these approaches seem to still be a long
way from resolving the weak uniqueness problem for the one-dimensional SPDE described in the
introduction which leads to much larger perturbtions.

\subsec{References}

\refer{ABP} S. Athreya, ÊR. F. Bass, and E.A. Perkins, H\"older norm estimates
for elliptic operators on finite and infinite dimensional spaces, 
{\sl Trans. Amer. Math. Soc.}, to appear.
Ê
\refer{B1} R.F. Bass, {\sl Probabilistic Techniques in Analysis.} Springer-Verlag, New York, 1995.

\refer{B2} R.F. Bass, {\sl Diffusions and Elliptic Operators.} Springer-Verlag,
New York, 1997.

\refer{BP} R.F. Bass and E. Perkins, Countable systems of degenerate stochastic differential 
equations with applications to super-Markov chains, {\sl Electron. J. Probab. \bf 9} (2004) 
634--673. 

\refer{Bi} P. Billingsley, {\sl Convergence of Probability Measures,
2nd ed.} Wiley, New York, 1999.

\refer{CD} P. Cannarsa and G. Da Prato, Infinite-dimensional elliptic equations with H\"older-continuous coefficients. {\sl Adv. Differential Equations \bf 1} (1996) 425--452.

\refer{Ce} S. Cerrai,
{\sl Second order PDE's in finite and infinite dimension.
A probabilistic approach.} 
Springer-Verlag, Berlin, 2001. 

\refer{D1} G. ÊDa Prato, Perturbation of Ornstein-Uhlenbeck semigroups. 
{\sl ÊRend. Istit. Mat. Univ. Trieste Ê\bf 28} Ê(1996) 101--126 (1997). 

\refer{D2} G. Da Prato, Bounded perturbations of Ornstein-Uhlenbeck semigroups. Ê
In: {\sl Evolution equations, semigroups and functional analysis (Milano, 2000)}, Ê97--114, 
Birk\-h\"aus\-er, Basel, 2002.

\refer{D3} G. Da Prato, Perturbations of Ornstein-Uhlenbeck operators: an analytic approach. Ê
In: {\sl Evolution equations: applications to physics, industry, life 
sciences and economics (Levico Terme, 2000)}, Ê127--139, ÊBirk\-h\"aus\-er, Basel, 2003.

\refer{D4} G. ÊDa Prato, ÊA new regularity result for Ornstein-Uhlenbeck generators and applications. 
{\sl ÊJ. Evol. Equ. Ê\bf Ê3} Ê(2003) 485--498.

\refer{DR} G. Da Prato and M. ÊR\"ockner, Singular dissipative stochastic equations in Hilbert spaces. Ê
{\sl Probab. Theory Rel. Fields Ê\bf 124} Ê(2002) 261--303. 

\refer{DL} G. Da Prato and A. Lunardi, On the Ornstein-Uhlenbeck operator in spaces of continuous functions. Ê
{\sl J. Funct. Anal. Ê\bf 131} Ê(1995) 94--114.

\refer{DZ} G. Da Prato and J. Zabczyk, {\sl Stochastic equations
in infinite dimensions.} Cambridge University Press,
Cambridge, 1992.

\refer{DM} D. Dawson and P. March, Resolvent estimates for Fleming-Viot operators and uniqueness of solutions to related martingale problems. {\sl J. Funct. Anal. } {\bf 132} (1995) 417-442.

\refer{KX} G. K. Kallianpur and J. Xiong, {\sl Stochastic Differential Equations in Infinite Dimensional Spaces.} IMS Lecture Notes-Monograph Series, Vol. 26,
IMS, Hayward, 1995.

\refer{L1} A. ÊLunardi, An interpolation method to characterize domains of generators of semigroups. Ê
{\sl Semigroup Forum Ê\bf 53} Ê(1996) 321--329. 

\refer{L2} A. Lunardi, Schauder estimates for a class of degenerate elliptic and parabolic 
operators with unbounded coefficients in $R\sp n$. Ê{\sl Ann. Scuola Norm. Sup. Pisa \bf Ê 24}
ÊÊ(1997) 133--164. 

\refer{L3} A. Lunardi, On the Ornstein-Uhlenbeck operator in $L\sp 2$ spaces with respect to 
invariant measures. Ê{\sl Trans. Amer. Math. Soc. Ê\bf 349} Ê(1997) 155--169. 

\refer {RN} F. Riesz, B. Sz.-Nagy, {\sl Functional Analysis}, Ungar, New York,
1955.

\refer{P} E.A. ÊPerkins, Dawson-Watanabe Superprocesses and Measure-Valued Diffusions, Ê 
In: {\sl Lectures on probability theory and statistics (St. Flour, 1999)}, 125--324,
Springer-Verlag, Berlin, 2002.

\refer{SV} D.W. Stroock and S.R.S. Varadhan, {\sl Multidimensional
Diffusion Processes}, Berlin, Springer, 1979.

\refer{Z} L. Zambotti, An analytic approach to existence and uniqueness for martingale problems in infinite dimensions. {\sl Probability Theory Êand Related Fields} {\bf 118} (2000) 147-168.

\bs

\noindent Addresses:

\smallskip\noindent S.R. Athreya:
Statmath Unit, Indian Statistical Institute, 8th Mile Mysore Road, R.V. College P.O.,
Bangalore 560059, India.

\smallskip\noindent R.F. Bass:
Department of Mathematics, University of Connecticut, \hfill\break
\hglue0.5in Storrs,
CT 06269, USA.

\smallskip\noindent M. Gordina:
Department of Mathematics, University of Connecticut, \hfill\break
\hglue0.5in Storrs,
CT 06269, USA.

\smallskip\noindent
E.A. Perkins:
Department of Mathematics, University of British Columbia,\hfill\break
\hglue0.5in Vancouver, B.C., Canada V6T 1Z2.

\bye